\documentclass{amsart}

\usepackage{amsmath}
\usepackage{amsfonts}
\usepackage{amscd}
\usepackage{amssymb}
\input xypic

\newtheorem{defn}{Definition}[section]
\newtheorem{thm}[defn]{Theorem}
\newtheorem{prop}[defn]{Proposition}
\newtheorem{lemma}[defn]{Lemma}
\newtheorem{eg}[defn]{Example}

\newcommand{\lm}{\ensuremath{\longrightarrow}}

\newcommand{\ra}{\rightarrow}

\DeclareMathOperator{\Hom}{\mbox{Hom}}

\DeclareMathOperator{\sspec}{\ensuremath{\mathcal{S}\mathit{pec}}}

\DeclareMathOperator{\Sym}{\mbox{Sym}}

\DeclareMathOperator{\Pic}{\mbox{Pic}\,}

\DeclareMathOperator{\e}{\varepsilon}

\DeclareMathOperator{\PP}{\mathbb{P}}
\DeclareMathOperator{\Z}{\mathbb{Z}}
\DeclareMathOperator{\cale}{\mathcal{E}}
\DeclareMathOperator{\call}{\mathcal{L}}

\DeclareMathOperator{\calo}{\mathcal{O}}

\DeclareMathOperator{\cala}{\mathcal{A}}

\DeclareMathOperator{\nez}{\overline{NE(Z)}}

\title{Numerically Calabi-Yau Orders on Surfaces}

\author{Daniel Chan}
\author{Rajesh S. Kulkarni}
\address{School of Mathematics, University of New South Wales, Sydney, NSW 2052, Australia}
\email{danielch@maths.unsw.edu.au}
\address{Department of Mathematics, Wells Hall, Michigan State University, East Lansing, MI 48824, USA}
\email{kulkarni@math.msu.edu}

\begin{document}
\begin{abstract}
This is part of an ongoing program to classify maximal orders on surfaces via their ramification data. Del Pezzo orders and ruled orders have been classified in \cite{CK, CIa} and \cite{AdJ}. In this paper, we classify numerically Calabi-Yau orders which are the noncommutative analogues of surfaces of Kodaira dimension zero.
\end{abstract}
\maketitle
Throughout, all objects and maps are assumed to
be defined over some algebraically closed base field $k$ of characteristic 0. 
\section{Introduction}
In \cite{A}, Artin raised the question of classifying noncommutative surfaces. This remains a difficult problem. Artin has initiated a program for classifying maximal orders on surfaces which is a far more tractable problem. There has been a lot of progress in this field over the past few years. Del Pezzo orders were classified in \cite{CK} and a type of minimal model program for orders was established in \cite{CIa} giving the fundamental dichotomy between orders which behave like Mori fibre spaces (that is, del Pezzo and ruled orders) and orders which are minimal models. Also, Artin and de Jong \cite{AdJ} have made substantial progress to the birational classification of orders on surfaces with results generalising Enriques' criterion for ruled surfaces. Ruled orders have been classified in terms of their ramification data (see \cite{CIa}). Maximal orders which behave like surfaces of Kodaira dimension zero is another class of interest and the aim of this note is to give a similar classification of these orders. 

More specifically, let $Z$ be a smooth projective surface and $X = \sspec \cala$ 
be a maximal order on $Z$. Let the ramification curves be $D_i$ and the corresponding ramification indices be 
$e_i$. It is convenient to define the canonical divisor of the order as 
\[K_X := K + \sum (1- \frac{1}{e_i})D_i  \]
where $K$ is the canonical divisor of the central surface (see \cite{CK} or \cite{CIa} for an explanation). It is a divisor on the central surface $Z$. We say that the order $X$ is {\em numerically Calabi-Yau} if $K_X$ is numerically trivial. 

By analogy with the commutative case, to classify noncommutative surfaces, we ought to impose some smoothness condition on the surface. The most useful such condition to date is that of a terminal order on a surface (also called stable in \cite{AdJ}). Following \cite{CIa}, we say that a maximal order $X$ as above is {\em terminal} if it satisfies the following conditions:
\begin{enumerate}
  \item Its centre $Z$ is smooth.
  \item The ramification divisor has only nodes as singularities.
  \item Suppose that $p$ is such a node, $D_i,D_j$ the ramification curves passing through 
        $p$ and $D_i',D_j'$ are the cyclic covers of $D_i,D_j$ describing the ramification. 
        Then either $e_i$ is the ramification index of $D_i'$ above $p$ or $e_j$ is the 
        ramification index of $D_j'$ above $p$.              
\end{enumerate}
Note that the definition depends purely on the ramification data of the order. The interested reader is encouraged to look at \cite{CIa} to see why this is a good class of maximal orders to look at. We will only mention here that, in a sense which can be made precise, every maximal order on a surface has a terminal resolution \cite{CIa}, Corollary~3.7. Hence the classification of maximal orders reduces to the classification of terminal orders the same way that the classification of surfaces reduces to that of smooth surfaces in the commutative case. 

Although our goal in this paper is to classify terminal orders which are numerically Calabi-Yau, we will need occasionally to consider the slightly more general class of canonical orders. Their definition, also motivated by Mori theory can be found in \cite{CIb}, \S2. We will only mention here that the definition depends only on the ramification data and the possible ramification data have been classified in \cite{CIb}, Theorem~3.5. Also, the centre of a canonical order has at worst canonical singularities. 

The key tool in the analysis is an exact sequence due to Artin and Mumford (see \cite{AM}, Theorem~1) which describes the Brauer group of the function field of a simply connected projective surface. It can be interpreted as a statement giving necessary and sufficient conditions for the existence of maximal orders with given ramification data (see \cite{CK}, Corollary~20 for a simple description of this). 

Our classification of the ramification data which occur for terminal numerically Calabi-Yau orders is as complete as one can hope for. The possible centres of such orders are 
\begin{enumerate}
  \item The quadric surface $\PP^1 \times \PP^1$, the Hirzebruch surface 
        $F_2 = \PP(\calo \oplus \calo(-2))$ or the blow-up of $\PP^2$ at up to 9 points in 
        almost general position.
  \item A surface $Z = \PP(\cale)$ ruled over an elliptic curve $C$ where either i) 
        $\cale \simeq \calo \oplus \call$ where $\call ^n \simeq \calo$ for 
        $n \in \{1,2,3,4\}$ or ii) $\cale$ is indecomposable of degree 1. 
  \item A surface of Kodaira dimension zero. 
\end{enumerate}
In the last case, the order is unramified and we give a description of the possible ramification curves and their ramification indices in the other two cases. In case~ii), there are a small number of possibilities for ramification data and so we are inclined to think that these orders are somehow special. 

{\it Conventions:} All cohomology groups are \'etale. The symbols $\equiv$ and $\sim$ denote numerical and linear equivalences respectively.

\section{Classification of Centres}
We showed in \cite{CK}, Theorem~12 that a del Pezzo order has 
del Pezzo centre. The key was that $-K_X$ positive implies that 
$-K$ is too. The result essentially then follows from the Nakai 
criterion. The proof works at the boundary too to show
\begin{thm}   \label{tcentre}  
If the anti-canonical divisor $-K_X$ of a maximal order is 
nef then so is $-K$. In particular, a numerically Calabi-Yau order on a 
surface is either an Azumaya algebra on a surface of 
Kodaira dimension zero or a ramified order on a ruled surface. 
\end{thm}
The following proposition allows one to reduce to the case 
where $Z$ is minimal.
\begin{prop}   \label{predtomin} 
Let $X'$ be a canonical numerically Calabi-Yau order on a smooth surface $Z'$. 
Then the terminal resolution is a numerically Calabi-Yau order. 
Conversely, suppose $X= \sspec \cala$ is a terminal 
numerically Calabi-Yau order on $Z$. Let $\pi: Z \lm Z'$ be a proper 
birational morphism to a surface $Z$ with rational singularities. Then 
$\pi_* \cala$ is a canonical numerically Calabi-Yau order on $Z'$. 
\end{prop}
\begin{proof} 
Let $X$ be the terminal resolution of a canonical numerically Calabi-Yau order and 
$Z$ be its centre. Write $\pi:Z \lm Z'$ for the natural map. Then 
$K_X = \pi^* K_{X'} \equiv 0$ so $X$ is numerically Calabi-Yau. Conversely, if $X$ 
is terminal numerically Calabi-Yau,  then $K_{X'} = \pi_* K_{X} \equiv 0$ since $Z$ has rational singularities (see \cite{KM}, Lemma~5.12). Also, for any exceptional curve 
$E \subset Z$ we have $K_X . E = 0$ so $X$ is canonical. 
\end{proof}
Let $X$ be a terminal numerically Calabi-Yau order. Suppose first that the centre $Z$ is rational. Theorem~\ref{tcentre} ensures that $-K$ is nef, so the proposition below and \cite{Dem}, Expos\'e III, \S 2 Th\'eor\`eme~1 imply that $Z$ is either i) the quadric surface $\PP^1 \times \PP^1$, ii) the Hirzebruch surface $F_2 = \PP(\calo \oplus \calo(-2))$ or iii) the blow-up of $\PP^2$ at up to 9 points in almost general position. 

We wish to determine which birationally ruled surfaces can be centres of terminal numerically Calabi-Yau orders. Let $Z'$ be a minimal model of $Z$ so that $Z'$ is a projective bundle of the form $\PP(\cale)$ where $\cale$ is a rank two 
vector bundle on a smooth curve $C$ of genus $g$. As in 
\cite{Hart}; Chapter~V we normalise $\cale$ so that $H^0(Z',\cale) \neq 0$ 
but $H^0(Z',\cale \otimes L) = 0$ for any invertible sheaf $L$ of 
negative degree. Let $\e:= -\deg \cale$. The Neron-Severi group is 
then freely generated by a fibre $F$ and a special section $C_0$ 
which has self-intersection $C_0^2 = -\e$. 
\begin{prop}  \label{pfewg} 
If $X$ is a terminal numerically Calabi-Yau order on a surface $Z$ 
then $K^2 \geq 0$. In particular, if $Z'$ is the minimal model as above then 
$g = 0$ or $1$. Moreover, $Z$ is already geometrically ruled if 
$g=1$. Suppose that $Z' = \PP(\cale)$ and $\e$ is the invariant defined above. 
When $g=0$, we have $0 \leq \e \leq 2$ and 
when $g=1$ we have $ \e = 0$ or $-1$. 
\end{prop}
\begin{proof} 
Since $K_X \equiv 0$ and $-K$ is nef we have 
\[  K^2 = -(\sum (1- \frac{1}{e_i})D_i ).K \geq 0  .\]
We also know for a geometrically ruled surface that 
$K^2 = 8(1-g)$. For this to be non-negative we must 
have $g = 0$ or 1. Now $K^2$ decreases by 1 on blowing up a 
closed point so we always have $g \leq 1$. 
Recall from \cite{Hart}, Chapter~V, Lemma~2.10 that 
\[  K \equiv -2C_0 + (2g - 2 - \e)F    .\]
Intersecting with $C_0$ and using $-K$ nef shows 
\[ 0 \leq -\e + 2 -2g    .\]
This and \cite{Hart}, Chapter~V Theorems~2.12 and 2.15 gives the bound on $\e$. 
\end{proof}
Given  a prime $p$ dividing some ramification index $e_i$, let 
$p^{\max}$ be the largest power of $p$ dividing any of the $e_i$. 
Let $C^p$ denote the union of ramification curves $D_i$ with 
$p^{\max} | e_i$. We recall a result of Artin's in 
\cite{CK}; Lemma~23.
\begin{lemma}  \label{lpa}  
Suppose $X$ is a maximal order on a smooth rational surface. Then 
$p_a(C^p) \geq 1$. 
\end{lemma} 
It is this corollary of \cite{AM}, Theorem~1 that we will mainly use to restrict possible ramification data. 

In the following, we will often record ramification data via ramification vectors, 
that is, we list $(e_1, e_2, \ldots )$ with 
$e_1 \leq e_2 \leq \ldots $ and the $e_i$'s are repeated with 
multiplicity. The multiplicity will depend on the case at hand, 
on $Z=\PP^2$ it will just be the degree of $D_i$. In general, the multiplicity will be of the form $D_i.H$ for some divisor $H$ in $Z$. There is a 
partial ordering on these vectors given by the product of the orders 
on $\Z$. Note that the function 
\[  \sum (1 - \frac{1}{e_i})   \]
is increasing with respect to this order. 
\section{Case I: Rational Centre}
\begin{prop}   \label{pplane0} 
A canonical numerically Calabi-Yau order on $\PP^2$ has ramification vector 
$(4,4,4,4)$ or $(2,6,6,6)$ or $(2,2,2,2,2,2)$. 
\end{prop} 
\begin{proof}
Let $d$ be the degree of the total 
ramification curve and write the ramification indices with 
multiplicity (equal to its degree for each component). Since $X$ is numerically Calabi-Yau we have 
\begin{equation}   \label{eplane}  
\Sigma:= \sum (1 - \frac{1}{e_i}) = 3  
\end{equation}
We see that $4 \leq d \leq 6$ and that the only solution 
when $d=6$ is $(2,2,2,2,2,2)$. If $d=5$ then as 
$e_1,e_2 \geq 2$ we see that 
\[ (1 - \frac{1}{e_3}) + (1 - \frac{1}{e_4}) + 
    (1 - \frac{1}{e_5}) \leq 2 \]
Hence $(e_3,e_4,e_5)$ is either a platonic triple 
$(2,2,n),(2,3,3)$, $(2,3,4),(2,3,5)$ or one of the triples 
$(3,3,3),(2,4,4),(2,3,6)$. We must have then $e_1=e_2 = 2$ 
from which we see that Lemma~\ref{lpa} is violated. Hence, there 
are no numerically Calabi-Yau orders on $\PP^2$ ramified on a quintic. 

Suppose now that $d = 4$. If $e_2 \geq 6$ then the only 
possibility is $(2,6,6,6)$. Hence by Lemma~\ref{lpa}, the only 
possible prime powers dividing the $e_i$ are $2,3,4,5$. If 
there are at least two primes dividing the $e_i$ then either 
the largest prime powers dividing $e_i$'s are 2,3, in which 
case we can only get $(2,6,6,6)$ or, the ramification vector is 
bounded below by either $(2,5,10,10)$ or $(3,4,12,12)$. Both 
these have $\Sigma > 3$. 
\end{proof}
From now on, we assume that $X$ is a terminal numerically Calabi-Yau order on a minimal surface $Z$. The next case we will examine is $Z = \PP^1 \times \PP^1$. We shall write ramification vectors with multiplicity $D_i.F$ where $F$ is one of the ruling fibres as specified in the case at hand. 
\begin{prop}   \label{pplane1}  
The terminal numerically Calabi-Yau orders on $\PP^1 \times \PP^1$ have ramification vectors $(3,3,3)$ for both the rulings  or $(2,4,4)$ or $(2,2,2,2)$ for any of the two natural rulings.
\end{prop}
\begin{proof}
Let $A$ and $B$ be the numerical classes of two rulings on $Z$. Then recall that the numerical class of the canonical divisor is
\[
K \equiv - 2 A - 2 B.
\]
Suppose we have a terminal numerically Calabi-Yau order on $Z$ ramified along irreducible curves $D_i$ with their respective ramification indices $e_i$. Let the numerical class of the curve $D_i$ be $a_i A + b_i B$. Then since
\[
K + \sum (1 - \frac{1}{e_i}) D_i \equiv 0,
\]
we get equations
\[
\sum (1 - \frac{1}{e_i}) a_i = 2 , \qquad \sum (1 - \frac{1}{e_i}) b_i = 2
\]
In particular, $\sum a_i \leq 4$ and $\sum b_i \leq 4$. Also clearly $\sum a_i \geq 3$ and $\sum b_i \geq 3$. Next we claim that the highest power $p^l$ of any prime $p$ dividing one of the $e_i$'s must be less than or equal to $4$. First we apply Lemma~\ref{lpa} in our situation. Let $p$ be a prime dividing one of the $e_i$'s. Let $C^p$ be the union of the ramification curves whose ramification indices are divisible by the highest power of $p$ dividing any of the $e_i$'s. Write $C^p \equiv \alpha A + \beta B$. Then by Lemma~\ref{lpa}, we get $p_a(C^p)= (\alpha - 1)(\beta - 1) \geq 1$. So $\alpha \geq 2,\ \beta \geq 2$. Now suppose that the highest power $p^l$ is bigger than 4. Using arguments above, we see that 
\[
2 = \sum (1 - \frac{1}{e_i}) a_i \geq \frac{8}{5} + \ {\rm other\ terms}
\]
Since $8/5 + 1/2 > 2$, we see that this is impossible. So now consider the case $\sum a_i = 3$. Let $e_1,e_2,e_3$ be the ramification indices written with multiplicities given by the $a_i$'s. They satisfy
\[
\frac{1}{e_1} + \frac{1}{e_2} + \frac{1}{e_3} = 1.
\]
So a simple calculation shows that in this case possible ramification vectors for the $A$-ruling are $(3,3,3)$ and $(2,4,4)$. Finally consider $\sum a_i = 4$. Then their ramification indices satisfy
\[
\frac{1}{e_1} + \frac{1}{e_2} + \frac{1}{e_3} + \frac{1}{e_4} = 2.
\]
In this case, the only ramification vector is $(2,2,2,2)$. The arguments work for the $B$-ruling as well. It is also easy to see that for these vectors, there is a terminal numerically Calabi-Yau order on $Z$ with these vectors as ramification vectors. 
It is also possible to have the ramification vector $(2,4,4)$ for one ruling and $(2,2,2,2)$ for the other ruling.
\end{proof}
We now consider the Hirzebruch surface $Z = \PP_{\PP^1}(\calo \oplus \calo(-2))$. Let the base curve be denoted by $C$ and the fibre of the projection $Z \ra C$ be denoted by $F$. We shall write ramification vectors with multiplicities given by $D_i.F$. 
\begin{prop}   \label{pplane}
Consider the surface $Z = \PP_{\PP^1}(\calo \oplus \calo(-2))$. 
The terminal numerically Calabi-Yau orders on $Z$ have ramification vectors $(2,4,4),(3,3,3)$ or $(2,2,2,2)$. Further description is given in the proof.
\end{prop}
\begin{proof}
The canonical divisor of $Z$ is
\[
K \equiv -2 C_0 - 4 F.
\]
Suppose a numerically Calabi-Yau order on $Z$ ramifies on irreducible curves $D_i \equiv a_i C_0 + b_i F$. Then the condition for the order to be numerically Calabi-Yau means 
\[
\sum (1 - \frac{1}{e_i}) a_i = 2, \qquad \sum (1 - \frac{1}{e_i}) b_i = 4.
\]
This gives the conditions
\[
3 \leq \sum a_i \leq 4, \qquad 5 \leq \sum b_i \leq 8. 
\]
Suppose the highest power for a prime $p$ dividing one of the $e_i$'s is $p^l$. Let $C^p$ be the union of the ramification curves whose ramification indices are divisible by $p^l$. Write $C^p \equiv a C_0 + b F$. Then by the adjunction formula, $p_a(C^p) = -a^2 + ab - b + 1$. By Lemma~\ref{lpa}, $-a^2 + ab - b + 1 \geq 1$. That is, $(a - 1)(b - (a + 1)) \geq 1$, so $ a \geq 2, \ b \geq a + 2$. Using this we see that 
\begin{eqnarray*}
2 & = & \sum (1 - \frac{1}{e_i}) a_i \\
  & \geq & (1 - \frac{1}{p^l})(2) + \frac{1}{2} (\sum a_i - 2) 
\end{eqnarray*}
So $ a \leq \sum a_i \leq 2 + \frac{4}{p^l}$. Similarly we get $ b \leq \sum b_i\leq 4 + \frac{8}{p^l}$. So unless $p^l = 2$, $a \leq 3$ and $b \leq 6$. If $p^l = 2$, then all the ramification indices are equal to $2$. In this case, the total ramification divisor  $D = \sum D_i = 4 C_0 + 8 F$. The ramification vector in this case is $(2,2,2,2)$. It is easy to see that there exist terminal numerically Calabi-Yau orders on $Z$ with this ramification data.

So now we consider $p^l \neq 2$ for some power of prime dividing one of the $e_i$'s. The considerations above show that $\sum a_i = 3$ and $\sum b_i = 5$ or $6$. By using the equation for the $a_i$'s above, we see that the ramification indices written with multiplicity $a_i$  must be of the form $(2,4,4)$, $(3,3,3)$ or $(2,3,6)$. 

We now use the condition $p_a(C^p) \geq 1$ to determine the possibilities in each of these 3 cases. To simplify notation, we write $D^e$ for the sum of ramification divisors $D_i$ with $e_i = e$. 
\begin{enumerate}
\item For the ramification vector $(2,4,4)$, we must have $D^4 \equiv 2C_0 + 4F, D^2 \equiv C_0 + 2F$. There is a ramification divisor with this ramification data. This follows since $h^0(Z, {\calo}(C_0 + 2 F)) = 4$, so there are many choices for components of this divisor. Then a ramification diagram follows by \cite{AM}, Theorem~1 or \cite{CK}, Corollary~20. 
\item For the ramification vector $(3,3,3)$ we have $D^3 \equiv 3C_0 + 6F$ and as in the previous case, there are lots of terminal orders with this ramification. 
\item For the ramification vector $(2,3,6)$ we must have $C_0 + 2F \equiv D^6 \equiv D^3 \equiv D^2$. A maximal order with this ramification data cannot be terminal.  
\end{enumerate}
This finishes the proof.
\end{proof}
\begin{eg}  \label{eg:enonminimal} 
Minimal terminal numerically Calabi-Yau order on a non-minimal surface. 
\end{eg}
Consider two smooth cubic curves $C_1,C_2$ on $\PP^2$. Assume 
that they intersect in 9 distinct points. Let $Z$ be the 
blowing up of $\PP^2$ at these 9 points and let $D_1,D_2$ 
be the strict transforms of $C_1,C_2$. Given \'etale 2-fold 
covers of $D_1,D_2$ we obtain from the Artin-Mumford sequence 
a numerically Calabi-Yau order on $Z$ ramified on $D_1,D_2$ with ramification index 
two. Contracting any of the (-1)-curves yields a non-terminal order 
so it is minimal. 
\section{Case II: Centre ruled over elliptic curves}
Next we consider now the case where $Z$ is a geometrically ruled over 
an elliptic curve $C$. In this case, the Artin-Mumford sequence has non-trivial cohomology. However, as has been noted in \cite{ford, ingalls}, this secondary obstruction to orders is not too difficult to work with. We summarise what we need in 
\begin{lemma}  \label{lcolin}
The cohomology in the Artin-Mumford sequence is $H^3(Z,\mu)$ which is isomorphic to $H^1(C, \mu)$. 
Let $\varphi: C' \lm C$ be a degree $n$ isogeny of elliptic curves. Then the kernel of the map induced by the Artin-Mumford sequence \[ H^1(C', \mathbb{Z}/n \mathbb{Z}) \stackrel{\gamma}{\longrightarrow} H^3(Z, \mathbb{Z}/n \mathbb{Z}) \xrightarrow{\sim} H^1(C, \mathbb{Z}/n \mathbb{Z}) \] is isomorphic to the kernel of $\varphi$ under the isomorphism $C'_n \cong H^1(C', \mathbb{Z}/n \mathbb{Z})$. Here $C'_n$ and $C_n$ are $n$-torsion points on respective elliptic curves. 
\end{lemma}
\begin{proof} 
The first part is Theorem 5, \cite{ford}. To prove the second part, we first note that the map $\gamma$ in the sequence above is the Gysin map (loc.~cit.) and that it is Poincar\'e dual to the restriction map $H^1(Z, \mathbb{Z}/n \mathbb{Z}) \stackrel{\gamma}{\longrightarrow} H^1(C', \mathbb{Z}/n \mathbb{Z})$. The last statement is \cite{milne}, Chapter VI, Remark (11.6b). Recall that the restriction map $H^1(C, \mathbb{Z}/n \mathbb{Z}) \longrightarrow H^1(Z, \mathbb{Z}/n \mathbb{Z})$ is an isomorphism by the Leray spectral sequence. So we get that the composite restriction map $H^1(C, \mathbb{Z}/n \mathbb{Z}) \longrightarrow H^1(C', \mathbb{Z}/n \mathbb{Z})$ is dual to the map of the statement. But this map can be identified with the map $C_n \longrightarrow C'_n$ obtained from the dual isogeny $C \longrightarrow C'$. Since the dual of this map is $C'_n \longrightarrow C_n$ obtained from the isogeny of the statement, we are done.
\end{proof}
\noindent
\textbf{Remark:} By the proposition below, the only possible ramification curves $D = \cup D_i$ will be disjoint unions of \'etale covers of $C$. Hence the possible ramification data of maximal orders with these ramification curves will correspond to elements of the kernel of 
\[ \oplus_i H^1(D_i,\mu) \lm H^1(C,\mu). \]
For a more general result, see \cite{AdJ}, Chapter~2. 
\vspace{2mm}
\begin{prop}   \label{pramiscan}  
If $Z$ is geometrically ruled over an elliptic curve and $X$ is terminal numerically Calabi-Yau then every ramification curve $D_i$ is numerically a 
multiple of $K= -2C_0 -\e F$. Hence, the curves $D_i$ are elliptic and disjoint, and the composite map $D_i \lm Z \lm C$ is an \'etale cover. 
\end{prop}
\begin{proof}
Since the base curve is elliptic we have 
\[  0 = K^2 = -(\sum (1- \frac{1}{e_i})D_i ).K  .\]
Now $-K.D_i \geq 0$ so we must have equality. Since $Z$ is 
ruled, the Kleiman-Mori cone $\nez$ is 2-dimensional and $K=0$ 
give extremal rays which include $K$ and $D_i$. The genus formula now gives $p_a(D_i) = 1$ which forces the ramification curves to be elliptic and \'etale covers of $C$. Also, $D_i.D_j = K^2 = 0$, so they are disjoint. 
\end{proof}
Below we write ramification vectors with multiplicities given by $D_i.F$. 
\begin{prop}  \label{pfewl}
Let $X$ be a terminal numerically Calabi-Yau order on 
$Z = \PP(\cale)$ where $\cale = \calo \oplus \call$ is a 
rank two vector bundle on an elliptic curve $C$ normalised as 
above. Then $\call^a \simeq \calo_C$ for some  
$a \in \{1,2,3,4\}$. The possible ramification vectors are $(2,2,2,2)$, $(3,3,3), (2,4,4), (2,3,6)$ ( cf.~the elliptic orders of \cite{AdJ}).

Conversely, the following combinations of ramification vectors and 
line bundles $\call$ have lots of terminal numerically Calabi-Yau orders. 
\begin{enumerate}
  \item $(2,2,2,2)$ and any $\call$ as above
  \item $(3,3,3)$ and $\call^a \simeq \calo$ for 
           $a \in \{1,2,3\}$.
  \item $(2,4,4)$ and $\call^a \simeq \calo$ for 
           $a = 1$ or 2.
  \item $(2,3,6)$ and $\call \simeq \calo$.
\end{enumerate}
\end{prop}
\begin{proof}
Since $\cale$ splits $\e=0$ and $K \sim -2C_0$. 
Let $\rho: Z \lm C$ be the natural map. From 
Proposition~\ref{pramiscan}, the ramification 
curves have the form $D_i \sim aC_0 + \rho^* P$ 
for some $a >0$ and some degree 0 divisor 
$P$ on $C$. Now 
\[  \sum (1- \frac{1}{e_i})D_i \equiv -K \sim 2C_0 \]
forces $a \leq 4$ and 
\begin{equation}  \label{eabig}
 \sum D_i . F \geq 3  
\end{equation}
Hence the only possibilities for ramification vectors are as listed above. 
We compute 
\begin{eqnarray*}
H^0(Z, \calo(aC_0 + \rho^* P)) & = & 
          H^0(C, \rho_*\calo_{\PP}(a) \otimes \calo_C(P)) \\
   & = &  H^0(C, \Sym^a(\cale) \otimes \calo_C(P))        \\
   & = &  H^0(C, (\calo \oplus \call \oplus \ldots \oplus \call^a) 
                                               \otimes \calo_C(P))
\end{eqnarray*}
This is non-zero if and only if $\calo_C(-P) \simeq \call^i$ for some 
$i \in \{0,\ldots, a\}$. Moreover, if $a = 1$ and 
$\call \not \simeq \calo$ then we see there are only 
two possible curves in the above linear systems, $C_0$ which 
corresponds to $\calo$ and another curve say $C_1$ corresponding to 
$\call$. Furthermore, if there are no isomorphisms between the 
powers of $\call$ then the only divisors in the linear 
systems above are linear combinations of $C_0,C_1$. In this case, 
$C_0,C_1$ are the only possibilities for the ramification curves 
$D_i$ so we obtain a contradiction to Equation~\ref{eabig}.  
Consequently, $\call^a \simeq \calo$ for some 
$a \in \{1,2,3,4\}$. 

The converse follows from Lemma~\ref{lcolin} and the remark following it. 
\end{proof}
\begin{thm}  
Let $Z = \PP(\cale)$ where $\cale$ is the extension of $\calo$ by 
$\calo$. Then there are no terminal numerically Calabi-Yau orders on $Z$. 
\end{thm}
\begin{proof}
Let $\cale$ be the nontrivial extension defined by the sequence
\begin{equation}
0 \ra \calo_C \ra \cale \ra \calo_C \ra 0.
\label{eq:def of e}
\end{equation}
So $\Lambda^2 \cale = \calo_C$ and the degree of $\cale$ is zero. This means $\e = 0$ and $K \sim -2C_0$. Now let $\rho: Z \ra C$ be the usual defining map. Then the ramification curves have to be of the form $D_i \sim a_iC_0 + \rho^* E$ for some integer $a_i > 0$ and $E$ is a divisor of degree $0$ on $C$. By Proposition~\ref{pramiscan}, $D_i$ is an elliptic curve. Hence it will be sufficient to prove that for any unramified map $C' \stackrel{g}{\ra} C$ of elliptic curves and any lift 
\[
\xymatrix{
   &       Z \ar[d]^\rho  \\
  C' \ar[ur]^{\tilde{g}} \ar[r]^g  & C  \\ 
}
\label{diagram:lift1}
\]
such that the the numerical class of the image of $C'$ in $Z$ is $a C_0,\ a \geq 1$ the image of $C'$ is actually $C_0$. First suppose there is an unramified map $C' \stackrel{g}{\ra} C$ of degree $a >1$ with a lift as in diagram \ref{diagram:lift1} which is of degree 1. By \cite{Hart}, Chapter 2, Proposition 7.12, such a lift gives an invertible sheaf $\call$ together with a surjective map $g^* \cale \ra \call$ of $\calo_{C'}$-modules. So we check for the possibilities of such invertible sheaves with a surjection as above. Note that $\det g^*\cale = \wedge^2 g^* \cale = g^* \wedge^2 \cale = \calo_{C'}$ . Also note from the same proposition that the invertible sheaf $\call$ will be the pull-back of $\calo(1)$ under $\tilde{g}$. Now if the image of $C'$ under $\tilde{g}$ is of numerical class $a C_0$, then it follows that the degree of $\call$ is $aC_0^2 = 0$. This then gives a short exact sequence
\begin{equation}
0 \ra {\call}^{\vee} \ra g^*\cale \ra {\call} \ra 0.
\label{eq:sequence for L}
\end{equation}
So now to get line bundles with such possible surjections, we compute ${\rm Hom}(g^*\cale, {\call})$ for invertible sheaves ${\call}$ of degree 0 on $C'$. Note that we have a short exact sequence 
\[
0 \ra \calo_{C'} \ra g^*\cale \ra \calo_{C'} \ra 0.
\]
In particular, $H^0(C', g^*\cale) \neq 0$ and is 2 or 1 depending on whether or not $g^*\cale$ splits as a direct sum of two copies of $\calo_{C'}$. But we claim that $g^*\cale$ does not split as this direct sum. To see this first note that
\begin{equation*}
g_*g^*\cale = \cale \otimes \; g_*\calo_{C'} = \cale \otimes (\oplus_{\alpha \in \hat{G}} {\call}_\alpha), 
\end{equation*}
where $\hat{G}$ is the kernel of the dual isogeny $C \ra C'$ and ${\call}_\alpha$ are invertible sheaves in the kernel of the homomorphism ${\rm Pic}(C) \ra {\rm Pic}(C')$. This shows that $\cale$ is a direct summand of $g_*g^*\cale$ with mupltiplicity one. Now for any degree 0 invertible sheaf ${\call} \neq \calo_{C}$, $H^0(C, \cale \otimes {\call}) = 0$, this follows from the long exact sequence of cohomology associated with the short exact sequence in \ref{eq:def of e} tensored with ${\call}$. Thus $h^0(C, \cale \otimes (\oplus_{\alpha \in \hat{G}} {\call}_\alpha)) = 1$ as $\cale$ is indecomposable. So, in particular, if $g^*\cale \cong \calo_{C'} \oplus \calo_{C'}$, then $h^0(C, g_*g^*\cale) = 2$. This is a contradiction. This also gives that $h^0(C', g^*\cale) = 1$

First consider the case ${\call} = \calo_{C'}$. So we need to compute ${\rm Hom}(g^*\cale, \calo_{C'})$. Using Serre duality, 
\[
{\rm Hom}(g^*\cale, \calo_{C'}) =  {\rm Ext}^0(g^*\cale , \calo_{C'}) = (H^1(C', g^*\cale))'
\]
Now $h^0(C', g^*\cale) = 1$ and Riemann-Roch for vector bundles on curves gives 
\[
h^0(C', g^*\cale) - h^1(C', g^*\cale) = 0
\] 
since the degree of $g^* \cale = 0$. So ${\rm Hom}(g^*\cale, \calo_{C'})$ is of dimension 1. But we already know one lift of $C' \ra Z$, namely, $C' \ra C \ra C_0$, where the last map is the section of $\rho$. So in such a case there are no emeddings of $C'$ in $Z$. 

Next consider the case when ${\call} \neq \calo_{C'}$. Using the long exact sequence on cohomology obtained from Eq.~\ref{eq:sequence for L}, we see that a short exact sequence such as in this equation cannot exist since $h^0(C', g^*\cale) \neq 0$. 

Finally we are left with the case when $a = 1$. By taking homomorphisms of the non-split exact sequence~(\ref{eq:def of e}) into $\call$ we see that there are no surjections $\cale \lm \call$ unless $\call \simeq \calo$ in which case we have $\Hom(\cale, \calo) = k$. Hence the only possible section of $Z \lm C$ is $C_0$. We conclude finally that the only possible ramification curve is $C_0$ which violates  the equality 
\[
\sum (1 - \frac{1}{e_i}) D_i \equiv 2 C_0.
\]
This completes the proof of the theorem. 
\end{proof}
\begin{thm}  
Let $X$ be a terminal numerically Calabi-Yau order on 
$Z = \PP(\cale)$ where $\cale$ is the non-split extension of a degree one line bundle $\call$ by $\calo$. Then the ramification indices are all equal to 2 and there is either i) a single ramification curve $D \equiv -2K$ or ii) two ramification curves $D_1,D_2$ with $D_i \equiv -K$. Conversely, there exist terminal numerically Calabi-Yau orders in both cases.
\end{thm}
\begin{proof} 
The strategy for the proof will be the same as in the last theorem. For an invertible sheaf ${\call}$ of degree 1 on $C$, let 
\[
0 \ra \calo_C \ra \cale \ra {\call} \ra 0
\]
be a non-split exact sequence of $\calo_C$-modules. Let $\rho: Z \ra C$ be the associated morphism. From the sequence, $\wedge^2 \cale = {\call}$ and so the degree of $\cale$ is 1 and $\e = - \deg \cale = -1$. So $K \equiv -2C_0 + F$. From Proposition \ref{pramiscan}, the ramification curves $D_i$ are of the form $D_i = a_i(-2C_0 + F) + \rho^* E$, where $E$ is a divisor of degree 0 on $C$. If $e_i$'s are the corresponding ramification indices, then we have the equality
\[
\sum (1- \frac{1}{e_i})D_i \equiv -K \equiv 2C_0 - F .
\]
Now substituting the numerical class of $D_i$ and then intersecting both sides with $F$, we get
\[
\sum (-a_i)(1 - \frac{1}{e_i}) = 1.
\]
Now since $-a_i \geq 1$ (by \cite{Hart}, Chapter V, Prop.~2.21) and $e_i \geq 2$, we get $\sum (-a_i) \leq 2$. In particular, the ramification divisor has the following possibilities:
\begin{enumerate}
\item There are two ramification curves $D_1, \; D_2$, their numerical class is $2 C_0 - F$ and their ramification indices are $e_1 = e_2 = 2$.
\item There is only one ramification curve $D_1$, its numerical class is $D_1 \equiv 4 C_0 - 2 F$ and the ramification index is $e_1 = 2$. 
\end{enumerate}
We will show that in both the cases, there exist ramification curves which yield terminal numerically Calabi-Yau orders as in the statement. We start by considering a degree 2 isogeny $g: C' \ra C$ which we will use to concoct a curve with numerical class $2C_0 - F$. Let $q$ be a closed point of $C'$ and $L = \calo_{C'}(q)$. First note that for any line bundle $F$ on $C$ of negative degree, we have $h^0(C,g_* L \otimes F) = h^0(C',L \otimes g^*F) = 0$ so $g_* L$ is normalised. Also, for any degree zero line bundle $F$ on $C$ we have $h^0(C,g_* L \otimes F) \neq 0$ so $g_* L$ is a non-split extension of a degree one line bundle $\call$ by $\calo$. In fact, if $N$ is the 2-torsion line bundle on $C$ defining the cover $C' \lm C$, then \cite{Hart}, Chapter~IV, Exercise~2.6 shows that 
\[ 
\call = \det g_* L = (\det g_* \calo) \otimes \calo(g(q)) = N \otimes \calo(g(q)). 
\]
Now a rank-$2$ vector bundle on $C$ that is an extension of a degree $1$ line bundle by $\calo$ is determined by its determinant. Hence, by choosing $q$ suitably, we may as well assume that $\cale = g_* L$. 

Let $L' = \calo(q')$ be the degree one line bundle on $C'$ which completes the canonical map $g^* g_* L \rightarrow L \rightarrow 0$ to an exact sequence below.
\begin{equation} \label{enonsplit}  
0 \lm L' \lm g^*g_* L \lm L \lm 0 .
\end{equation}
Computing determinants one finds 
\[ L' \otimes L \simeq g^* \det g_*L \simeq g^* N \otimes g^*\calo(g(q)) \simeq 
\calo(g^{-1}(g(q))). \]
Hence $q, \ q'$ are in the same fibre of $g$. Since $L,\ L'$ are non-isomorphic, $g^* \cale \simeq L \oplus L'$ and there are precisely two lifts of $g$ to maps $\tilde{g}:C' \lm Z$. They correspond to the surjections $\cale \lm L$ and $\cale \lm L'$. Next since the degree of $L$ is 1, the corresponding lift has degree 1 as well. This follows since if the degree of this lift is $d$, and the degree of the restriction of $\calo(1)$ to $\tilde{g}(C')$ is $d'$, then the degree of $L$ is $dd'$. 

Finally, let $a C_0 + b F$ be the numerical class of the image of such a lift. Since $g$ is of degree 2, it follows that $a=F.(a C_0 + b F) = 2$. Since the line bundle $L$ is of degree 1 and is a pull-back of $\calo(1)$, $C_0.(a C_0 + b F) = 1$. That is, $a + b = 1$, and hence $b = -1$. Now we have three curves $C_i$ with degree $2$ isogenies $g_i: C_i \ra C$. So we get embeddings associated with these curves, $\tilde{g}_i: C_i \ra Z$. We claim that the images $\tilde{g}_i(C_i)$ are distinct. Note that there are three line bundles $L_i$ on $C$ associated with these covers such that $g_j^* L_i$ is trivial if and only if $i = j$. Now $\rho \tilde{g}_i = g_i$, so $g_i^* L_i = \tilde{g}_i^* \rho^* L_i$, and $g_i^* L_i$ is trivial if and only if $(\rho|_{\tilde{g}_i(C_i)})^* L_i$ is trivial. So if any two of the images $\tilde{g}_i(C_i)$ are the same, the pull-back of the corresponding $L_i, \ L_j$ are trivial under $g_i$ which is a contradiction.

Now by Lemma~\ref{lcolin} and the remark following it, there are choices for double covers of a pair of such curves which lead to maximal orders. 

We now construct maximal orders with ramification as in case ii). With the above notation, let $f:C'' \lm C'$ be the degree 2 isogeny which annihilates the 2-torsion point $q-q'$ of $\Pic C'$. We let $h = gf$ and pull-back equation~(\ref{enonsplit}) by $f$ to obtain $h^*\cale \simeq f^*L \oplus f^* L'$. By our choice of $f$, $f^*L \simeq f^*L'$, and so ${\rm Hom}(h^*\cale, f^*L)$ is of dimension 2. But ${\rm Hom}(g^*\cale, L)$ is of dimension 1 as is seen by applying ${\rm Hom}(-, L)$ to Eq.~\ref{enonsplit}. So there must exist some embedding of $C''$ into $Z$. Again it is easy to see that the numerical class of the image of such an embedding must be $4 C_0 - 2 F$. Any smooth divisor in such a class gives rise to terminal numerically numerically Calabi-Yau orders by Lemma~\ref{lcolin} and the subsequent remark.
\end{proof}

\end{document}